\documentclass{amsart}
\usepackage[english]{babel}
\usepackage[cp1250]{inputenc}
\usepackage{amsmath,amssymb}
\usepackage{microtype}
\usepackage{float}
\usepackage{enumitem}
\usepackage{mathrsfs}

\usepackage{tikz}
\usetikzlibrary{matrix,arrows}
\usetikzlibrary{decorations.markings}

\newcommand\blfootnote[1]{%
  \begingroup
  \renewcommand\thefootnote{}\footnote{#1}%
  \addtocounter{footnote}{-1}%
  \endgroup
}

\newtheorem{thm}{Theorem}[section]
\newtheorem{prop}[thm]{Proposition}
\newtheorem{lem}[thm]{Lemma}

\newtheorem*{claim*}{Claim}

\theoremstyle{definition}
\newtheorem{definition}[thm]{Definition}
\newtheorem{remark}[thm]{Remark}

\newtheorem{notation}[thm]{Notation}

\newtheorem*{remark*}{Remark}

\def\ZZ{\mathbb Z}

\def\QQ{\mathbb Q}

\DeclareMathOperator{\tb}{tb}
\DeclareMathOperator{\rot}{rot}

\DeclareMathOperator{\AbMap}{AbMap}
\DeclareMathOperator{\Map}{Map}

\definecolor{b}{RGB}{27,77,137}
\definecolor{o}{RGB}{211,84,0}
\definecolor{g}{RGB}{13,113,3}

\begin{document}

\title[Fillability of small Seifert fibered spaces]{Fillability of small Seifert fibered spaces}

\author{Irena Matkovi\v{c}}
%\address{Alfr\'ed R\'enyi Institute of Mathematics, Budapest, Hungary}
\address{Mathematical Institute, University of Oxford, Oxford, UK}
\email{irena.matkovic@maths.ox.ac.uk}

\begin{abstract}
On small Seifert fibered spaces $M(e_0;r_1,r_2,r_3)$ with $e_0\neq-1,-2,$ all tight contact structures are Stein fillable. This is not the case for $e_0=-1$ or $-2$. However, for negative twisting structures it is expected that they are all symplectically fillable. Here, we characterize fillable structures among zero-twisting contact structures on small Seifert fibered spaces of the form $M(-1;r_1,r_2,r_3)$. The result is obtained by analyzing monodromy factorizations of associated planar open books.
\end{abstract}

%\subjclass[2010]{57R17}

\maketitle

%%%%%%%%%%%%%%%%%%%%%%%%%%%%%%%%%%%%%%%%%%%%%%%%%%%%%%%%%%%%%%

\section{Introduction}%%%%%%%%%%%%%%%%%%%%%%%%%%%%%%%%%%%%%%%%%%%%%%%%%%%%

Seifert fibered $3$-manifolds not carrying fillable contact structures have been singled out by Lecuona and Lisca \cite{LL}; they call them manifolds of special type. In this paper, we are interested in exactly which contact structures on small Seifert fibered spaces are fillable. For the surgery presentation of the underlying manifold see Figure \ref{fig:e_0}.

\begin{figure}[h]%----------------------------------------------------------------------------
\begin{tikzpicture}
\draw[white] (7.5,-0.5) circle (0pt) node[black,left,scale=0.9]{$e_0$};
\draw (9,0) ellipse (1.8cm and 0.6cm);
\draw [fill=white,white] (8.35,-0.45) rectangle (8.55,-0.65);
\draw [fill=white,white] (9.3,-0.6) rectangle (9.45,-0.4);
\draw [fill=white,white] (10.1,-0.35) rectangle (10.3,-0.55);
\draw (7.9,-.55) arc (-160:170:0.3);
\draw[white] (8,-0.8) circle (0pt) node[below,black,scale=0.9]{$-\frac{1}{r_1}$};
\draw (8.8,-.7) arc (-160:160:0.3);
\draw[white] (8.9,-.9) circle (0pt) node[below,black,scale=0.9]{$-\frac{1}{r_2}$};
\draw (9.7,-.65) arc (-160:170:0.3);
\draw[white] (9.85,-.9) circle (0pt) node[below,black,scale=0.9]{$-\frac{1}{r_3}$};
\end{tikzpicture}
\caption{Small Sefert fibered space $M(e_0;r_1,r_2,r_3)$.}
\label{fig:e_0}
\end{figure}%--------------------------------------------------------------------------------

Tight contact structures on small Seifert fibered spaces $M(e_0;r_1,r_2,r_3)$ where $e_0\in\ZZ$ and $r_i\in\QQ\cap(0,1),$ are completely classified \cite{Wu, GLS0} whenever $e_0\neq -1$ or $-2$; they are all given by Legendrian surgery on the Stein fillable $S^3$, hence they are also Stein fillable. The same holds for the Seifert manifolds with $e_0=-2$ which are $L$-spaces \cite{G}. On $M(-1;r_1,r_2,r_3)$ there are essentially two types of tight contact structures, distinguished by the maximal twisting of the regular fiber; that is, the maximal difference between the contact framing and the fibration framing within the smooth isotopy class of the fiber. The negative twisting structures are related to the transverse contact structures, and they are expected to be all at least symplectically fillable. On the other hand, the zero-twisting tight contact structures share a common contact surgery description \cite{LS.III} and are conjecturally \cite{S} characterized by non-vanishing of the Ozsv\'ath-Szab\'o contact invariant $c(M,\xi)\in\widehat{HF}(-M,\mathbf t_\xi)$; in the particular case of $L$-spaces this covers all tight structures and it has been confirmed in \cite{M}.

\begin{figure}[h]%----------------------------------------------------------------------------
\begin{tikzpicture}
\begin{scope}[scale=1.5, decoration={markings,mark=at position 0.5 with {\arrow{>}}}]
\begin{scope}[shift={(0,0)}]
\draw [domain=-0.5*pi:1.5*pi, scale=0.5] plot (\x, {2*sin(\x r)});
\draw [scale=0.5,postaction={decorate}] (-0.5*pi,-2)--(1.5*pi,-2);
\draw [scale=0.5] (-0.5*pi,-2) circle (0.1pt) node[left,scale=0.9]{$-\frac{1}{r_3}$};
\end{scope}
\begin{scope}[shift={(.06,-.3)},scale=0.95]
\draw [domain=-0.5*pi:1.5*pi, scale=0.5] plot (\x, {2*sin(\x r)});
\draw [scale=0.5,postaction={decorate}] (-0.5*pi,-2)--(1.5*pi,-2) node[right,scale=0.9]{$-\frac{1}{r_2}$};
\end{scope}
\begin{scope}[shift={(.12,-.6)},scale=0.9]
\draw [domain=-0.5*pi:1.5*pi, scale=0.5] plot (\x, {2*sin(\x r)});
\draw [scale=0.5,postaction={decorate}] (-0.5*pi,-2)--(1.5*pi,-2);
\draw [scale=0.5] (-0.5*pi,-2) circle (0.1pt) node[left,scale=0.9]{$-\frac{1}{r_1}$};
\end{scope}
\begin{scope}[shift={(.18,-.9)},scale=0.85]
\draw [domain=-0.5*pi:1.5*pi, scale=0.5] plot (\x, {2*sin(\x r)});
\draw [scale=0.5,postaction={decorate}] (-0.5*pi,-2)--(1.5*pi,-2)node[right,scale=0.9]{$+1$};
\end{scope}
\begin{scope}[shift={(.24,-1.2)},scale=0.8]
\draw [domain=-0.5*pi:1.5*pi, scale=0.5] plot (\x, {2*sin(\x r)});
\draw [scale=0.5,postaction={decorate}] (-0.5*pi,-2)--(1.5*pi,-2)node[right,scale=0.9]{$+1$};
\end{scope}
\end{scope}
\end{tikzpicture}
\caption{Contact structures on $M(-1;r_1,r_2,r_3)$.}
\label{fig:SFS}
\end{figure}%---------------------------------------------------------------------------------

Zero-twisting tight contact structures on $M(-1;r_1,r_2,r_3)$ are all described by contact surgery diagrams of Figure \ref{fig:SFS}, as shown by Lisca and Stipsicz in \cite{LS.III}. Recall that such a diagram gives a family of contact structures, whose elements can be specified by replacing each contact $-\frac{1}{r_i}$-surgery with a Legendrian surgery along a chain $L_i$ -- called a leg -- of unknots $v_i^j$ whose Thurston-Bennequin invariants are determined by the continued fraction expansion of \[-\frac{1}{r_i}=-a^0_i-\frac{1}{\ddots-\frac{1}{-a^{k_i}_i}}=[a^0_i,\ldots,a^{k_i}_i],\ a^j_i\geq 2,\]
as \[\begin{array}{l}\tb_i^0=-a_i^0 \text{ for the leading unknot of } L_i \text{ and}\\ \tb_i^j=-a_i^j+1 \text{ for } j>0, \end{array}\]
and rotation numbers are chosen arbitrarily in \[\rot_i^j\in\{\tb_i^j+1,\tb_i^j+3,\dots,-\tb_i^j-1\}.\]

These structures are all supported by planar open books (see Subsection \ref{Ss2.1}); but in contrast to contact structures on small Seifert spaces with $e_0\neq -1$, not all tight ones are Stein fillable. With the aid of a theorem of Wendl \cite{W} (see also Theorem \ref{thm:W}), we see that non-Stein fillable structures are not fillable at all. Non-fillability was first observed by Ghiggini, Lisca and Stipsicz in \cite{GLS} for a particular structure on $M(-1;\frac{1}{2},\frac{1}{2},\frac{1}{p})$; moreover, they obtained the classification of tight structures on $M(-1;r_1,r_2,r_3)$ for $r_1\geq r_2\geq\frac{1}{2}$. Based on their classification, Plamenevskaya and Van Horn-Morris \cite{PVH-M} recognized exactly which of these Seifert manifolds admit non-fillable tight structures; this was achieved by using Wendl's work and obstructing existence of positive factorization for the abelianized monodromy of the standardly associated (planar) open books. On the other hand, Lecuona and Lisca \cite{LL} showed that, when $M(-1;r_1,r_2,r_3)$ is an $L$-space and $r_i+r_j<1$ for all pairs $i,j$, topology (the diagonalization argument) prevents existence of Stein fillings.

Here, we show that all fillable zero-twisting structures on $M(-1;r_1,r_2,r_3)$ arise as Legendrian surgeries on the tight $S^1\times S^2$. For $L$-spaces this covers all fillable structures, and hence implies the result of Lecuona and Lisca.
More specifically, we show that fillability of a given surgery presentation is completely decided on specific sublinks representing $S^1\times S^2$, whose tightness is in turn met by a unique choice of rotation numbers for this sublink.

First notice that, whenever $r_i+r_j\geq 1$ there exists a  truncated continued fraction $-\frac{1}{s_i}=[a^0_i,\dots,a^{m_i}_i]<[a^0_i,\dots,a^{k_i}_i]=-\frac{1}{r_i}$ with $m_i\leq k_i$, and for $r_j$ alike, such that $s_i+s_j=1$ (see \cite[Lemma 3.2]{LL}). We will call any subsurgery of the contact surgery presentation in Figure \ref{fig:SFS} which consists of the two unstabilized unknots with $+1$-coefficient and two truncated legs, representing rational numbers $-\frac{1}{s_i}$ such that $s_i\leq r_i, s_j\leq r_j$ and $s_i+s_j=1$, a {\em circular sublink}.

Additionally, we will say that a Legendrian knot is {\em fully positive} if all its stabilizations are positive, that is $\rot=-(\tb+1)$. When all the knots forming a leg are fully positive, the leg will be said to be positive. Analogously, we define a {\em fully negative} Legendrian knot and a negative leg.
Now, a circular sublink whose one leg is positive and the other one negative will be refered to as a {\em balanced sublink}.
With the terminology set we can state our result.

\begin{thm}\label{thm}
Assume that a contact structure $\xi$ on $M(-1;r_1,r_2,r_3)$ is given by some surgery diagram of Figure \ref{fig:SFS}. 	Then $\xi$ is fillable if and only if the surgery presentation contains a balanced sublink.
\end{thm}

As explained above, the theorem covers all zero-twisting tight structures on $M(-1;r_1,r_2,r_3)$ and all tight contact structures when the underlying manifold is also an $L$-space. Additionally, when $r_3=0$ (equivalently, when there is no surgery along $L_3$) and $r_1+r_2=1$ we have the following.

\begin{prop}\label{prop:tightSxS}
Contact surgery diagram as in Figure \ref{fig:SFS} describes the tight contact $S^1\times S^2$ if and only if it equals some balanced link.
\end{prop}

\subsection*{Overview} In Section \ref{Sec2}, we recall how to associate open books to the given surgery presentations, and present main properties of planar monodromies. The proof of Theorem \ref{thm} is split between the following sections. In Section \ref{Sec3}, we show that the surgery along any balanced link indeed gives the fillable $S^1\times S^2$. In Section \ref{Sec4}, we obtain negative results by obstructing positive factorization of monodromy in the abelianization of the mapping class group of the planar page.

\subsection*{Acknowledgement} I thank András Stipsicz and Alberto Cavallo for their interest and the valuable comments on an early draft of this paper. My research has been supported by the NKFIH Grant \'Elvonal (Frontier) KKP 126683 and by the European Research Council (ERC) under the European Union's Horizon 2020 research and innovation programme (grant agreement No 674978).

\section{Open book presentation}\label{Sec2}%%%%%%%%%%%%%%%%%%%%%%%%%%%%%%%%%%%

\subsection{Planar open book from the contact surgery presentation}\label{Ss2.1}
Recall that Legendrian surgeries of Figure \ref{fig:SFS} can be given on a planar page of the associated open book, describing its monodromy, as follows. Look at Figure \ref{fig:view} top for the corresponding illustration.

One $+1$-surgery along an unknot with $\tb=-1$ is presented by an annulus with identity monodromy, the other $+1$-surgery manifests in a negative Dehn twist along its core. 
\begin{notation}
Write $\pi^\text{in}$ and $\nu^\text{out}$ for the inner and the outer boundary of the annulus. Additionally, write $N$ for the curve which supports the negative Dehn twist and, abusing the notation, also for the negative Dehn twist itself.
\end{notation}

Every other unknot contributes a positive Dehn twist on a stabilized annulus. Concretely, we insert a hole, encircled by one boundary-parallel positive Dehn twist, for every stabilization of every unknot in the surgery diagram; the stabilization holes which correspond to positive stabilizations lie between the inner boundary $\pi^\text{in}$ of the annulus and its core $N$, the negative ones between $N$ and the outer boundary $\nu^\text{out}$. 
\begin{notation}
Denote by $\nu^j_i$ and $\pi^j_i$ any of the stabilization holes which correspond to negative and positive stabilizations, respectively, of the unknot $v_i^j$. 
When grouped into certain types, we use $\nu_i$ for any of $\cup_{j}\nu^j_i$, similarly $\nu^{\geq j}_i$ for any of $\cup_{y\geq j}\nu^y_i$, and $\nu$ to denote any of $\nu^\text{out}\cup\nu_i$; the notation for the $\pi$-type holes is analogous.
Note that, since the names of holes (equivalently, boundary components) are chosen as common names, we will refer to a single (not specified) hole with the given name $\chi$ as a $\chi$-hole.
\end{notation}
\begin{remark}
Using $|\cdot |$ for the number of respective holes, we see that $2+|\nu^j_i|+|\pi^j_i|$ equals $a^j_i$ for $j>0$ and $a^0_i+1$ for $j=0$, $-1-|\nu^j_i|-|\pi^j_i|=\tb^j_i$ and $|\pi^j_i|-|\nu^j_i|=\rot^j_i$. 
\end{remark}

From the leading unknot $v^0_i$ of each leg we get a positive Dehn twist along a push-off of the core $N$ modified by encircling an additional $\nu_i^0$-hole for each negative stabilization, and avoiding a $\pi_i^0$-hole for each positive stabilization. The twists corresponding to the subsequent unknots $v_i^j$ in each leg are obtained from a push-off of the twist corresponding to the preceding unknot $v_i^{j-1}$, additionally encircling all $\nu_i^j$-holes and avoiding all $\pi_i^j$-holes.
\begin{notation}
Denote $T_i^j$ the curve which corresponds to the unknot $v_i^j$. Abusing the notation, we will use $T_i^j$ also for the corresponding positive Dehn twist. As in the case of holes, we will need also common names such as $T_i=\cup_j T_i^j$, referring to any single one of them as a $T_i$-twist or curve.
\end{notation}
\begin{remark}
The unknots $T_i^j$ and $v_i^j$ are not exactly the same; while $v_i^j$ with fixed $i$ form a chain $L_i$, the corresponding $T_i^j$ give the rolled-up diagram of $L_i$ (as described in \cite{Sch}). Anyway, since the Legendrian push-off and the meridian of a Legendrian knot are Legendrian isotopic after we perform $-1$-surgery on the original knot \cite{DG}, the two presentations give the same contact manifold. 
\end{remark}

All together, the associated open book has as a page the punctured disk with one more hole than the number of stabilizations of all unknots in the surgery presentation, and  the monodromy (up to conjugation) is given as a product of the negative Dehn twist $N$, positive Dehn twists $T_i^j$ for all $j$ for all $i$ and a single boundary twist along all boundary components except $\pi^\text{in}$ and $\nu^\text{out}$.

\begin{figure}%--------------------------------------------------------------------
\begin{tikzpicture}
\begin{scope}[shift={(0,1)}, scale=1.55,decoration={markings, mark=at position 0.88 with {\arrow{>}}}]
\begin{scope}[gray!70]
\draw (0,0) circle (2.57cm);\fill (1.9,-1.86) circle (0pt) node[right, scale=.7]{$\nu^\text{out}$};\draw (0,0) circle (0.15cm) node[scale=.7]{$\pi^\text{in}$};
\draw (.7,.67) circle (0.15cm) node[scale=.7]{$\pi^0_1$}; \draw (-.65,.56) circle (0.15cm) node[scale=.7]{$\pi^0_2$}; \draw (.30,.35) circle (0.15cm) node[scale=.7]{$\pi^1_1$}; 
\draw (-.5,-1.41) circle (.15cm) node[scale=.7]{$\nu^0_3$};\draw (.5,-1.41) circle (.15cm) node[scale=.7]{$\nu^0_3$};\draw (1.7,-1) circle (.15cm) node[scale=.7]{$\nu^1_3$};
\draw (1.88,1.31) circle (.15cm) node[scale=.7]{$\nu^1_2$};
\end{scope}
\draw (.7,.67) circle (0.18cm); \draw (-.65,.56) circle (0.18cm); \draw (.30,.35) circle (0.18cm); 
\draw[dashed,postaction={decorate}] (0,.1) circle (1.3cm); \draw[thick,color={g},postaction={decorate}] (0,-.1) circle (1.71cm); \draw[thick,color={g},postaction={decorate}] (.07,-.2) circle (2.12cm); 
\draw[thick,color={o},rotate around={42:(.43,.33)},postaction={decorate}] (.43,.33) ellipse (.74cm and .4cm);
\draw[thick,color={o},rotate around={33:(.7,.7)},postaction={decorate}] (.77,.6) ellipse (1.5cm and .66cm);
\draw[thick,color={b},rotate around={140:(-.4,.33)},postaction={decorate}] (-.4,.33) ellipse (.7cm and .3cm);\draw[thick,color={b},rotate around={140:(-.4,.33)},postaction={decorate}] (-.4,.33) ellipse (.75cm and .35cm);
\draw[thick,color={b},rotate around={147:(-.17,.3)},postaction={decorate}] (-.11,.31) ellipse (.94cm and .61cm);
\draw (-.5,-1.41) circle (.18cm); \draw (.5,-1.41) circle (.18cm); \draw (1.7,-1) circle (.18cm); \draw (1.88,1.31) circle (.18cm);
\end{scope}

\begin{scope}[shift={(0,-7.5)}, scale=1.55]
\begin{scope}[gray!70]
\draw (.5,0) circle (2.59cm);\fill (-1.57,-1.95) circle (0pt) node[right, scale=.7]{$\nu^0_3$};\draw (0,0) circle (0.15cm) node[scale=.7]{$\pi^\text{in}$};
\draw (.7,.67) circle (0.15cm) node[scale=.7]{$\pi^0_1$}; \draw (-.65,.56) circle (0.15cm) node[scale=.7]{$\pi^0_2$}; \draw (.30,.35) circle (0.15cm) node[scale=.7]{$\pi^1_1$}; 
\draw (.5,-1.47) circle (.15cm) node[scale=.7]{$\nu^0_3$};\draw (1.93,-.7) circle (.15cm) node[scale=.7]{$\nu^1_3$};
\draw (1.91,1.34) circle (.15cm) node[scale=.7]{$\nu^1_2$};
\draw (2.5,0) circle (.17cm) node[scale=.7]{$\nu^\text{out}$};
\end{scope}
\draw (.5,0) circle (2.55cm);
\draw (.7,.67) circle (0.18cm); \draw (-.65,.56) circle (0.18cm); \draw (.30,.35) circle (0.18cm); 
\draw[dashed] (0,.1) circle (1.3cm); 
\draw[thick,color={g},rotate around={115:(2.2,.7)}] (2.2,.7) ellipse (1cm and .4cm);
\draw[thick,color={g},rotate around={95:(2.1,.33)}] (2.1,.37) ellipse (1.39cm and .73cm);
\draw[thick,color={o},rotate around={42:(.43,.33)}] (.43,.33) ellipse (.74cm and .4cm);
\draw[thick,color={o},rotate around={33:(.7,.7)}] (.77,.6) ellipse (1.5cm and .66cm);
\draw[thick,color={b},rotate around={140:(-.4,.33)}] (-.4,.33) ellipse (.7cm and .3cm);
\draw[thick,color={b},rotate around={140:(-.4,.33)}] (-.4,.33) ellipse (.75cm and .35cm);
\draw[thick,color={b},rotate around={147:(-.17,.3)}] (-.11,.31) ellipse (.94cm and .61cm);
\draw (.5,-1.47) circle (.18cm); \draw (1.93,-.7) circle (.18cm); \draw (1.91,1.34) circle (.18cm);
\end{scope}
\end{tikzpicture}
\caption{Illustration of our notation conventions on an example: $\tb_1=(-2,-2,-1)\ \text{and}\ \rot_1=(1,1,0),\ \tb_2=(-2,-2)\ \text{and}\ \rot_2=(1,-1),\ \tb_3=(-3,-2)\ \text{and}\ \rot_3=(-2,-1).$ In gray are boundary components of the punctured disk. The full curves correspond to positive Dehn twists, the boundary twists of the stabilization holes are black, the $T_1$-twists are blue, the $T_2$-twists orange and the $T_3$-twists green. The dashed curve represents the negative Dehn twist $N$. 
The page is shown in two perspectives: with initial outer boundary (top) and with outer boundary in one $\nu^0_3$-hole (bottom).}
\label{fig:view}
\end{figure}%------------------------------------------------------------------

\subsection{Planar monodromy}
Since our contact structures are all planar, the following theorem of Wendl ensures that to prove non-fillability it suffices to study positive factorizations of the given monodromy.

\begin{thm}\cite[Corollary 2]{W}\label{thm:W}
A planar contact manifold is strongly symplectically (and thus Stein) fillable if and only if every supporting planar open book has monodromy isotopic to a product of positive Dehn twists.
\end{thm}

Let us briefly review the characteristic features of the abelianized planar mapping classes, as used by Plamenevskaya and Van Horn-Morris in \cite{PVH-M}.

The mapping class group of a planar surface (in the presentation of Margalit and McCammond \cite{MMc}) is described (geometrically) on a disk, $\mathbf D_n$, with $n$ holes arranged in the roots of unity. The group $\Map\mathbf D_n$ is generated by all convex Dehn twists (that is, the twists whose core is the boundary of the convex hull of a set of holes), and factored by commutators of disjoint twists and all lantern relations. Then, up to conjugation we have the following.
\begin{lem}
A Dehn twist as an element of $\AbMap\mathbf D_n$  is determined by the set of holes it encircles. 
\end{lem}

Furthermore, any mapping class $\phi$ factors into a product of Dehn twists, and each Dehn twist can be, using the lantern relations, decomposed into pairwise (around a pair of holes) and boundary (around a single hole) Dehn twists; when a positive Dehn twist encircles $r$ holes, it provides $r-1$ positive pairwise twists and $r-2$ negative boundary twists, both around each of its holes. Hence, if we define the {\em single and pairwise multiplicities}, $m_\alpha(\phi)$ and $m_{\alpha\beta}(\phi)$, to be the number of twists (counted with signs) on the disk with all but one hole $\alpha$, or a pair of holes $\alpha$ and $\beta$, capped off, these contain a complete homological information about $\phi$.
\begin{lem}\cite[page 2084]{PVH-M}
A mapping class $\phi$ as an element of $\AbMap\mathbf D_n$ is uniquely determined by a collection of multiplicities $\{m_\alpha(\phi),m_{\alpha\beta}(\phi)\}$.
\end{lem}

In particular, in a positive factorization, the number of non-boundary twists around every hole is bounded from above by the number (counted with signs) of all twists encircling this hole in any given presentation. 

In the following, we will make extensive use of an iterated lantern relation, also known as a {\em daisy relation}, which we state in the lemma below and illustrate in Figure \ref{fig:daisy}.

\begin{figure}[h]%---------------------------------------------------------------
\begin{tikzpicture}
\begin{scope}[scale=1.4]

\begin{scope}[shift={(-1,0)}]
\begin{scope}[gray!70]
\draw[thick] (0,0) circle (1.9cm);\fill (-1.3,0) circle (.07cm);
\fill (1.3,0) circle (.07cm); \fill (1.2,.5) circle (.07cm); \fill (.9,1) circle (.07cm); \fill (.9,-1) circle (.07cm); \fill (1.2,-.5) circle (0pt) node[scale=.9]{\reflectbox{$\ddots$}};
\end{scope}
\fill[white] (-1,-1.4) circle (0pt) node[black,right,scale=.7]{$B$};
\draw (-1.3,0) circle (.1cm); \draw (-1.3,0) circle (.14cm); \draw (-1.3,0) circle (.23cm) node[below,scale=.7]{$\ \ \ \ \ \ \ \ \ \ \ B_0$}; \fill (-1.12,.07) circle (0pt) node[scale=.7]{$\vdots$};
\draw (1.3,0) circle (.1cm) node[right,scale=.7]{$\ \ B_3$}; \draw (1.2,.5) circle (.1cm) node[right,scale=.7]{$\ \ B_2$}; \draw (.9,1) circle (.1cm) node[right,scale=.7]{$\ \ B_1$}; \draw (.9,-1) circle (.1cm) node[right,scale=.7]{$\ B_{k+1}$};
\draw (0,0) circle (1.83cm);
\end{scope}

\fill (1,0) circle (0pt) node[right]{$\leftrightarrow$};
\begin{scope}[shift={(3.4,0)}]
\begin{scope}[gray!70]
\draw[thick] (0,0) circle (1.9cm);\fill (-1.3,0) circle (.07cm);
\fill (1.3,0) circle (.07cm) node[black,below,scale=.7]{$\ \ \ \ \ \ \  A_3$}; \fill (1.2,.5) circle (.07cm) node[black,below,scale=.7]{$\ \ \ \ \ \ \  A_2$}; \fill (.9,1) circle (.07cm) node[black,below,scale=.7]{$\ \ \ \ \ \ \  A_1$}; \fill (.9,-1) circle (.07cm) node[black,above,scale=.7]{$\ \ \ \ \ \ \ \ \ \ \ A_{k+1}$}; \fill (1.2,-.5) circle (0pt) node[scale=.9]{\reflectbox{$\ddots$}};
\end{scope}
\draw[rotate around={24.5:(-1.3,0)}] (-.1,0) ellipse (1.4cm and .3cm);
\draw[rotate around={11.5:(-1.3,0)}] (0,0) ellipse (1.45cm and .3cm);
\draw (0,0) ellipse (1.5cm and .3cm);
\fill (.8,-.43) circle (0pt) node[scale=.9]{$\vdots$};
\draw[rotate around={-24.5:(-1.3,0)}] (-.1,0) ellipse (1.4cm and .3cm);
\draw (1,0) ellipse (.63cm and 1.3cm); \fill[white] (.7,-1.4) circle (0pt) node[black,right,scale=.7]{$C$};
\end{scope}

\end{scope}
\end{tikzpicture}
\caption{Daisy relation.}
\label{fig:daisy}
\end{figure}%----------------------------------------------------------------------

\begin{lem}\cite[Lemma 3.5]{PVH-M}\label{lem:daisy}
In the mapping class group of the disk with $k+2$ holes, the positive Dehn twists $B, B_0,B_1,\dots, B_{k+1}$ and $A_1,\dots,A_{k+1},C$, as denoted in Figure \ref{fig:daisy}, satisfy the relation
\[(B_0)^kB_1\cdots B_{k+1}B = CA_{k+1}\cdots A_1\]
\end{lem}

\begin{remark}
The daisy relation of the disk with $k+2$ holes exactly describes the rational blow-down along $L((k+1)^2,k)$, as monodromy substitution for the Lefschetz fibration \cite{EMVH-M}.
\end{remark}

\section{Surgery links of tight $S^1\times S^2$}\label{Sec3}%%%%%%%%%%%%%%%%%%%%%%%%%%%%%%%%%

\begin{lem}
The contact surgery presentation given by a circular link smoothly describes $S^1\times S^2$.
\end{lem}

\proof
The circular link  smoothly consists of four $-1$-linked unknots with framing coefficients $0,0,-\frac{s_1+1}{s_1},-\frac{s_2+1}{s_2}$ for some $s_1+s_2=1$. Switching to integral coefficients, we have the legs $L_1$ and $L_2$ in place of rational framed unknots, with surgery coefficients $(-a_i^0-1,-a_i^1,\dots,-a_i^{m_i})$ where $-\frac{1}{s_i}=[a^0_i,\dots,a^{m_i}_i]$ for $i=1,2$. Blowing-up once at the linking point (followed by a blow-down of the two $(+1)$-framed meridians of the thus-added curve), we obtain a chain of unknots with coefficients $(-a^{m_1}_1,\dots,-a^0_1,-1,-a^0_2,\dots,-a^{m_2}_2)$. Since $[a^{m_1}_1,\dots,a^0_1,1,a^0_2,\dots,a^{m_2}_2]=0$, this chain can be successively,  starting with the middle $-1$-surgery, blown-down ending in a $0$-framed unknot.
\endproof

\begin{remark}\label{rmk:dual}
Notice that, after the blow-up the two legs of a circular link become dual to each other (that is, they describe a lens space and its orientation reversal). Explicitly, the coefficients of the two are related as follows (here, $-2^{\times b}$ means a chain of $b$-many unknots with framing $-2$): 
\[\begin{array}{lllcccc}L'_1&:&(-b_1-2,&-2^{\times b_2},& -b_3-3,&\dots,&-2^{\times b_m}\text{ or }-b_m-2)\\ L'_2&:&(-2^{\times b_1},&-b_2-3,&-2^{\times b_3},&\dots,& -b_m-2\text{ or }-2^{\times b_m})\end{array}.\]

Hence, since the legs are affected by the blow-up only at the leading unknots, the smooth surgery coefficients of the original legs correlate as:
\[\begin{array}{lllcccc}L_1&:&(-b_1-3,&-2^{\times b_2},& -b_3-3,&\dots,&-2^{\times b_m}\text{ or }-b_m-2)\\ L_2&:&(-3,-2^{\times (b_1-1)},&-b_2-3,&-2^{\times b_3},&\dots,& -b_m-2\text{ or }-2^{\times b_m})\end{array}.\]
\end{remark}

\begin{prop}\label{prop:+factor}
The contact surgery presentation given by a balanced link corresponds to the tight $S^1\times S^2$.
\end{prop}

\proof
We prove that the presented contact manifold is Stein fillable by describing a concrete positive factorization of the associated monodromy.

Since a balanced link is circular, we can write out the smooth surgery coefficients of its two legs as in Remark \ref{rmk:dual}:
\[\begin{array}{lllcccc}L_1&:&(-b_1-3,&-2^{\times b_2},& -b_3-3,&\dots,&-2^{\times b_m}\text{ or }-b_m-2)\\ L_2&:&(-3,-2^{\times (b_1-1)},&-b_2-3,&-2^{\times b_3},&\dots,& -b_m-2\text{ or }-2^{\times b_m})\end{array}\]
for some $b_i\geq 0$.
Without loss of generality, we choose $L_1$ to be negative and $L_2$ positive. 

Recall from Section \ref{Ss2.1} that the associated monodromy factorizes into a product of the negative Dehn twist $N$, the positive Dehn twist $T_i^j$ for every unknot  $v_i^j$ and the positive boundary twists of stabilization holes. In the case of a balanced link with $L_1$ negative and $L_2$ positive, we have only $\nu_1$- and $\pi_2$-stabilization holes, all $T_1$-curves lie outside $N$,  while all $T_2$-curves lie inside. 
We can rewrite this monodromy by repeated use of the daisy relation as follows; look also at the example given by Figure \ref{fig:+factor}.

\begin{figure}%--------------------------------------------------------------------
\begin{tikzpicture}
\begin{scope}[scale=1.5]
\begin{scope}[shift={(-1,0)}]
\begin{scope}[gray!70]
\draw[thick] (0,0) circle (1.9cm);\fill (0,0) circle (.07cm);
\fill (.5,.5) circle (.07cm); \fill (.5,-.3) circle (.07cm); \fill (-.5,-1.25) circle (.07cm); \fill (1.3,0) circle (.07cm); \fill (1.6,-.5) circle (.07cm);
\end{scope}
\draw[color={o},dashed] (0,0) circle (1.1cm); \draw (0,0) circle (1.53cm); \draw (0,0) circle (1.85cm); 
\draw[color={o},rotate around={-35:(.23,-.13)}] (.23,-.13) ellipse (.5cm and .3cm);
\draw (0,0) circle (.1cm); \draw (0,0) circle (.14cm);
\draw[color={o}] (.5,.5) circle (.1cm); \draw (.5,-.3) circle (.1cm); \draw (-.5,-1.25) circle (.1cm); \draw[color={o}] (1.3,0) circle (.1cm); \draw (1.6,-.5) circle (.1cm);
\end{scope}

\fill[white] (1.2,0) circle (0pt) node[black,right]{$\leftrightarrow$};
\begin{scope}[shift={(3.3,0)}]
\begin{scope}[gray!70]
\fill (1.8,-.7) circle (.07cm);\fill (0,0) circle (.07cm);
\fill (.5,.5) circle (.07cm); \fill (.5,-.3) circle (.07cm); \fill (-.5,-1.25) circle (.07cm); \fill (1.3,0) circle (.07cm); \fill (1.6,-.5) circle (.07cm);
\end{scope}
\draw[color={o},dashed] (0.6,0) circle (.9cm); \draw (0,0) circle (1.53cm); \draw (1.8,-.7) circle (.1cm); 
\draw[color={o},rotate around={-30:(.73,.43)}] (.91,.37) ellipse (.65cm and .3cm);
\draw[color={o}] (.6,-.05) ellipse (.8cm and .4cm);
\draw (0,0) circle (.1cm); \draw (0,0) circle (.14cm);
\draw (.5,-.3) circle (.1cm); \draw (-.5,-1.25) circle (.1cm); \draw (1.6,-.5) circle (.1cm);
\end{scope}

\fill[white] (3.3,-2) circle (0pt) node[black,right]{$\|$};
\begin{scope}[shift={(3.3,-4)}]
\begin{scope}[gray!70]
\fill (1.8,-.7) circle (.07cm);\fill (0,0) circle (.07cm);
\fill (.5,.5) circle (.07cm); \fill (.5,-.3) circle (.07cm); \fill (-.5,-1.25) circle (.07cm); \fill (1.3,0) circle (.07cm); \fill (1.6,-.5) circle (.07cm);
\end{scope}
\draw[dashed] (0.6,0) circle (.9cm); \draw[color={o}] (0,0) circle (1.53cm); \draw (1.8,-.7) circle (.1cm); 
\draw[color={o},rotate around={-30:(.73,.43)}] (.91,.37) ellipse (.65cm and .3cm);
\draw (.6,-.05) ellipse (.8cm and .4cm);
\draw (0,0) circle (.1cm); \draw (0,0) circle (.14cm);
\draw[color={o}] (.5,-.3) circle (.1cm); \draw (-.5,-1.25) circle (.1cm); \draw (1.6,-.5) circle (.1cm);
\end{scope}

\fill[white] (1.1,-4) circle (0pt) node[black,right]{$\leftrightarrow$};
\begin{scope}[shift={(-1.7,-4)}]
\begin{scope}[gray!70]
\fill (1.8,-.7) circle (.07cm);\fill (0,0) circle (.07cm);
\fill (.7,.44) circle (.07cm); \fill (.7,-.3) circle (.07cm); \fill (-.1,-1.25) circle (.07cm); \fill (1.3,0) circle (.07cm); \fill (1.6,-.5) circle (.07cm);
\end{scope}
\draw[color={o}] (.88,.1) circle (.53cm); \draw[color={o},rotate around={-45:(1.2,-.1)}] (1.27,-.1) ellipse (.95cm and .25cm);
\draw[color={o},rotate around={-20:(1.2,-.4)}] (1.2,-.47) ellipse (.85cm and .35cm);
\draw[dashed] (0.6,0) circle (.9cm); \draw[color={o},dashed] (1.2,-.2) circle (1cm);  \draw (1.8,-.7) circle (.1cm); 
\draw (.6,-.05) ellipse (.8cm and .4cm);
\draw (0,0) circle (.1cm); \draw (0,0) circle (.14cm);
\draw (-.1,-1.25) circle (.1cm); \draw (1.6,-.5) circle (.1cm);
\end{scope}

\fill[white] (-1.3,-6) circle (0pt) node[black,right]{$\|$};
\begin{scope}[shift={(-1.7,-8)}]
\begin{scope}[gray!70]
\fill (1.8,-.7) circle (.07cm);\fill (0,0) circle (.07cm);
\fill (.7,.44) circle (.07cm); \fill (.7,-.3) circle (.07cm); \fill (-.1,-1.25) circle (.07cm); \fill (1.3,0) circle (.07cm); \fill (1.6,-.5) circle (.07cm);
\end{scope}
\draw[color={o}] (.88,.1) circle (.53cm); \draw[rotate around={-45:(1.2,-.1)}] (1.27,-.1) ellipse (.95cm and .25cm);
\draw[rotate around={-20:(1.2,-.4)}] (1.2,-.47) ellipse (.85cm and .35cm);
\draw[color={o},dashed] (0.6,0) circle (.9cm); \draw[color={o},dashed] (1.2,-.2) circle (1cm);  \draw[color={o}] (1.8,-.7) circle (.1cm); 
\draw (.6,-.05) ellipse (.8cm and .4cm);
\draw[color={o}] (0,0) circle (.1cm); \draw[color={o}] (0,0) circle (.14cm);
\draw[color={o}] (-.1,-1.25) circle (.1cm); \draw[color={o}] (1.6,-.5) circle (.1cm);
\end{scope}

\fill[white] (1.1,-8) circle (0pt) node[black,right]{$\leftrightarrow$};
\begin{scope}[shift={(2.9,-8)}]
\begin{scope}[gray!70]
\fill (1.8,-.7) circle (.07cm);\fill (0,0) circle (.07cm);
\fill (.7,.44) circle (.07cm); \fill (.7,-.3) circle (.07cm); \fill (1.3,0) circle (.07cm); \fill (1.5,-.5) circle (.07cm);
\draw[thick] (.7,0) circle (1.57cm);
\end{scope} 
\draw[rotate around={-45:(1.2,-.1)}] (1.27,-.1) ellipse (.95cm and .25cm);
\draw[rotate around={-20:(1.2,-.4)}] (1.2,-.47) ellipse (.85cm and .35cm);
\draw (.6,-.05) ellipse (.8cm and .4cm);
\draw[color={o}] plot [smooth cycle] coordinates {(-.11,.05)(.5,.7)(.9,.9)(1.3,.9)(1.5,.5)(1.9,-.7)(1.7,-.7)(1.3,.7)(.9,.8)(.5,.6)(.11,-.05)};
\draw[color={o}] plot [smooth cycle] coordinates {(-.13,.03)(.5,-.7)(1.1,-.9)(1.62,-.5)(1.4,-.43)(1.1,-.7)(.5,-.6)(.11,.03)};
\end{scope}
\end{scope}
\end{tikzpicture}
\caption{Example of positive factorization: $\tb_1=(-3,-2)$ with $\rot_1=(-2,-1)$ and $\tb_2=(-2,-2,-1)$ with $\rot_2=(1,1,0)$. On the first and the last picture the page is presented as a punctured disk with outer boundary in $\nu^\text{out}$ and one of $\nu^0_1$, respectively. Intermediate steps are presented as punctured spheres. In each row, the twists involved in a single application of the daisy relation are highlighted in orange.}
\label{fig:+factor}
\end{figure}%------------------------------------------------------------------

For the ease of notation, we write $v(b_\ell)$ for the unknot with the surgery coefficient $-b_\ell-3$ for $\ell<m$ and $-b_m-2$ for $\ell=m$. So, in our general notation $v(b_\ell)$ equals $v_1^{(\sum_{l=1}^{\ell'}b_{2l})+\ell'}$ for odd $\ell=2\ell'+1$ and $v_2^{(\sum_{l=1}^{\ell'}b_{2l-1})+\ell'-1}$ for even $\ell=2\ell'$. We attune the notation for twists and stabilization holes so that, the twist $T(b_\ell)$ corresponds to the unknot $v(b_\ell)$, and the holes $\nu(b_\ell)$ or $\pi(b_\ell)$ correspond to its stabilizations.

Throughout, we imagine the page as a disk with the outer boundary in one of the $\nu(b_1)$-holes; this hole (and its boundary parallel twist)  will not be considered a stabilization hole and hence, we change its notation to $\delta$. Then, all the $T_1$-curves encircle $\nu^\text{out}$, and all the $T_2$-curves encircle $\pi^\text{in}$. The new twists which arise by applying lantern relations, will be described by the subset of holes they encircle on the disk bounded by $\delta$.

To obtain a positive factorization, we will need $m$ applications of the daisy relation. It will be alternately applied from inside, involving some $T_2$-twists, and from outside, involving some $T_1$-twists.

For the zeroth application of the daisy relation (from inside), we consider:
\begin{itemize}
\item the first $b_1$ parallel $T_2$-twists;
\item the boundary twists of $b_1$ of the $\nu(b_1)$-holes (note that, the only non-considered $\nu(b_1)$-hole we set as the outer boundary);
\item the boundary twist of the $\pi_2^0$-hole;
\item the negative Dehn twist $N$.
\end{itemize}
The daisy relation (Lemma \ref{lem:daisy}), for which the $T_2$-twists take the role of $B_0$, the stabilization holes the role of $B_1,\dots,B_{b_1+1}$ and $N$ the role of $(A_{b_1+1})^{-1}$, results in:
\begin{itemize}
\item a new negative twist $N_0$ around $\{\pi^\text{in}\cup\pi_2\cup\nu(b_1)\}$, playing the role of $(B)^{-1}$, while we eliminate the negative twist $N$;
\item a new positive twist $D_0$ around all considered stabilization holes $\{\pi_2^0\cup\nu(b_1)\}$, playing the role of $C$;
\item new positive twists in role of $A_1,\dots,A_{b_1}$, which we will not keep track of because they remain unchanged in the continuation.
\end{itemize}

For the first application of the daisy relation (from outside), we consider:
\begin{itemize}
\item the positive Dehn twist $T(b_1)=T_1^0$ and all its parallel $T_1$-twists; all together there are $b_2+1$ of them;
\item the boundary twists of all $b_2+1$ of $\pi(b_2)$-holes;
\item the positive Dehn twist $D_0$.
\end{itemize}
We apply the daisy relation on them, so that the $T_1$-twists take the role of $B_0$ and the stabilization holes together with $D_0$ the role of $B_1,\dots,B_{b_2+2}$. This results in:
\begin{itemize} 
\item a new negative twist $N'_1$ around $\{\nu^\text{out}\cup\nu_1\cup\pi_2^0\cup\pi(b_2)\}$, playing the role of $(B)^{-1}$;
\item a new positive twist $D_1$ around $\{\pi_2^0\cup\nu(b_1)\cup\pi(b_2)\}$, playing the role of $C$, while we eliminate the positive twist $D_0$;
\item new positive twists in role of $A_1,\dots,A_{b_2+2}$, which we will again not keep track of.
\end{itemize}

We continue by alternately applying the daisy relation from inside and from outside. The $\ell^\text{th}$ application involves $T(b_\ell)$ and all its parallel twists, along with the stabilization holes of the unknot $v(b_{\ell+1})$. From inside (for the even applications $\ell=2\ell'$), $T(b_\ell)$ and its parallels are $T_2$-twists and the stabilization holes are $\nu(b_{\ell+1})$-holes; the daisy relation affects also $N_{\ell-2}$ and $D_{\ell-1}$ which get cancelled and replaced by enlarged curves $N_\ell$ and $D_\ell$, additionally encircling $\nu(b_{\ell+1})$-holes. From outside (for the odd applications $\ell=2\ell'+1$), $T(b_\ell)$ and its parallels are $T_1$-twists and the stabilization holes are $\pi(b_{\ell+1})$-holes; the daisy relation affects also $N'_{\ell-2}$ and $D_{\ell-1}$ which get cancelled and replaced by enlarged curves $N'_\ell$ and $D_\ell$, additionally encircling $\pi(b_{\ell+1})$-holes. 
So, after the $\ell^\text{th}$ application of the daisy relation, the twists contain:
\[\begin{array}{ll}
 \ell=2\ell':     &   D_\ell  =\{\pi^0_2\cup\nu(b_1)\cup\pi(b_2)\cup\cdots\cup\pi(b_{\ell})\cup\nu(b_{\ell+1})\}  \\
   & N_\ell  =\{\pi^\text{in}\cup\pi_2\cup\nu(b_1)\cup\dots\cup\nu(b_{\ell+1})\} \\
  & \\
   \ell=2\ell'+1:     &   D_\ell  =\{\pi^0_2\cup\nu(b_1)\cup\pi(b_2)\cup\cdots\cup\nu(b_{\ell})\cup\pi(b_{\ell+1})\}  \\
     & N'_\ell  =\{\nu^\text{out}\cup\nu_1\cup\pi_2^0\cup\pi(b_2)\cup\dots\cup\pi(b_{\ell+1})\} \\
\end{array}\]

Note that, as in the first two applications which we have explicitly described above, the $T$-twists always take the role of $B_0$, the stabilization holes together with $D_{\ell-1}$ the role of $B_1,\dots,B_{b_{\ell+1}+2}$, and $N_{\ell-2}$ or $N'_{\ell-2}$ the role of $(A_{b_{\ell+1}+2})^{-1}$; while for the resulting twists, $D_\ell$ takes the role of $C$ and $N_\ell$ or $N'_\ell$ the role of $(B)^{-1}$. Note also that, $N_\ell$ exists only for even $\ell$ and $N'_\ell$ only for odd $\ell$, either of them remaining untouched by the $(\ell+1)^\text{th}$ application of the daisy relation.

Finally, in the last, the $(m-1)^\text{th}$, application of the daisy relation, there are $b_m+1$ parallel twists $T(b_{m-1})$, but there are only $b_m$ stabilization holes of $v(b_m)$, so together with $D_{m-2}$ only $b_m+1$ twists in role of $B_i$ for $ i>0$. Hence, we involve as an additional $B_i$, the twist $T(b_m)$, which in our perspective appears as the boundary twist around $\nu^\text{out}$ when $m$ odd and around $\pi^\text{in}$ when $m$ even. So, after we apply the daisy relation, the twist $D_{m-1}$ encircles:
\[\begin{array}{ll}
 D_{m-1}  =\{\pi^0_2\cup\nu(b_1)\cup\pi(b_2)\cup\cdots\cup\pi(b_{m-1})\cup\nu(b_{m})\cup\nu^\text{out}\} & \text{for odd } m  \\
  & \\
    D_{m-1}  =\{\pi^0_2\cup\nu(b_1)\cup\pi(b_2)\cup\cdots\cup\nu(b_{m-1})\cup\pi(b_{m})\cup\pi^\text{in}\} & \text{for even } m \\
\end{array}\]
It exactly agrees with $N'_{m-2}$ for odd $m$ and $N_{m-2}$ for even $m$, and hence the corresponding negative twist gets cancelled by $D_{m-1}$. The other negative twist, $N_{m-1}$ when $m$ odd and $N'_{m-1}$ when $m$ even, encircles all the holes and it cancels with the positive Dehn twist along the outer boundary $\delta$.
\endproof

\section{Obstructing positive factorizations}\label{Sec4}%%%%%%%%%%%%%%%%%%%%%%%%%%%%%%%%%%%

When there is no balanced sublink of the surgery presentation, we show that the contact manifold of Figure \ref{fig:SFS} cannot be fillable.

\begin{prop}\label{prop:no+factor}
When the contact surgery presentation of Figure \ref{fig:SFS} does not contain a balanced sublink, the corresponding monodromy does not admit any positive factorization. 
\end{prop}

The rest of the section is devoted to the proof of Proposition \ref{prop:no+factor}. First, we recall that the oppositely stabilized leading unknots are needed already for tightness (Proposition \ref{prop:opposite}). Then, after setting up some further conventions and notation, we undertake a systematic analysis of possible positive factorizations of the monodromy $\phi$ as being read from the surgery presentation, obtaining eventually that the factorization cannot exist in the absence of a balanced sublink. First, we study how possible positive factorizations behave with respect to the $\pi$-holes (Lemmas \ref{l1.1}--\ref{l1.2}). Then, we specify particular properties of the twists which encircle $\nu$-holes (Lemmas \ref{l2}--\ref{pr2}). Finally, we obtain that (the lifts of) the factorizations as characterized among the $\pi$-holes can fulfill the above list of properties only when there is a balanced sublink (Lemmas \ref{l3.1}--\ref{l4}).

\vspace{.2cm}\noindent\emph{Necessary condition for tightness.}\vspace{.1cm}

We can argue via convex surface theory that some surgery configurations are always overtwisted.

\begin{prop}\label{prop:opposite}
Necessarily for tightness, the presentation admits a leg with fully negative leading unknot and a leg  with fully positive leading unknot.
\end{prop}

\proof
The absence of a pair of fully oppositely stabilized leading unknots is a very special case of the overtwistedness conditions given in \cite[Section 5]{M}.
\endproof

Due to Proposition \ref{prop:opposite}, we may, up to an overall orientation of the surgery link, assume that there is only one leg, say $L_3$, whose leading unknot is fully negative.

\vspace{.2cm}\noindent\emph{Conventions and notation.}\vspace{.1cm}

Throughout the proof, we try to build a positive factorization of the monodromy $\phi$ only on the level of abelianization. Abusing the notation, all mapping classes and their factorizations are considered as elements of $\AbMap$. In particular, we are interested in Dehn twists only up to conjugation, and we freely choose the order of Dehn twist factors. Note that we reserve the capital Greek letters to denote specific factorizations, in contrast to the maps as a whole.

It will be essential to look at a planar page of the open book from different {\em perspectives}, by which we mean a diffeomorphism of the page where a fixed boundary component becomes the outer boundary of the disk. 
Our preferred perspective will be the punctured disk obtained by setting one of the $\nu^0_3$-holes to be the outer boundary; call it $\mathbf D=\mathbf D_{\pi\cup\nu}$. As already in the proof of Proposition \ref{prop:+factor}, this $\nu^0_3$-hole (and its boundary parallel twist) will not be considered a stabilization hole and hence, we change its notation to $\delta$. When some of the holes in $\mathbf D$ are capped off, we denote this by putting the remaining holes in the index; for example, the notation $\mathbf D_\chi$ means the page $\mathbf D$ with all but the $\chi$-holes capped off. In arguments, we will interchangeably use two other perspectives on the (possibly capped-off) page: the {\em initial} with $\nu^\text{out}$ as the outer boundary of the disk, and the {\em turned-over} with $\pi^\text{in}$ as the outer boundary. Same (collections of) twists will be denoted by the same names regardless the perspective.

The multiplicities with respect to each perspective will be denoted by capital $M$ in $\mathbf D$, by $m$ for the initial disk, and by $m'$ for the turned-over one. When we wish to emphasize from which factorization the multiplicity was read, we put the factorization in the parenthesis, for example $M(\Psi)$. 

Let $\Phi$ denote the original factorization of the monodromy $\phi$ as being read from the surgery presentation. We recall from Section \ref{Ss2.1} that $\Phi=N\cdot\prod_{i,j}T_i^j\cdot D_\partial$ where $N$ is the negative twist representing one $+1$-surgery, the positive twists $T_i^j$ correspond to $-1$-surgeries on $v_i^j$, and $D_\partial$ is a boundary twist about all but two, $\nu^\text{out}$ and $\pi^\text{in}$, boundary components.

\vspace{.2cm}\noindent\emph{Positive factorizations and $\pi$-holes.}\vspace{.1cm}

To begin with, let us study how positive factorizations can possibly behave with respect to the $\pi$-holes.

\begin{lem}\label{l1.1} 
By capping off all the $\nu$-holes, we descend from $\AbMap \mathbf D_{\pi\cup\nu}$ to $\AbMap \mathbf D_{\pi}$, sending $\phi$ to $\overline{\phi}$. This maps the given factorization $\Phi$ to $\overline\Phi$, which is a composition $\overline{\Phi}_1\overline{\Phi}_2\overline{\Phi}_3$ with $\overline{\Phi}_i$ being a product of the  $T_i$-twists and the boundary twists around the $\pi_i$-holes. Every positive factorization $\overline\Psi$ of  $\overline\phi$ splits into subfactorizations $\overline\Psi=\overline{\Psi}_1\overline{\Psi}_2\overline{\Psi}_3$ so that $\overline{\Psi}_i$ and $\overline{\Phi}_i$ describe the same element $\overline{\phi}_i$ in $\AbMap\mathbf D_{\pi}$.
\end{lem}

\proof
The $\overline\Phi$ itself presents a positive factorization of the restricted monodromy $\overline\phi$. Indeed, the only negative twist of $\Phi$ cancels with the boundary twist of $\delta$ after we have capped-off the $\nu$-holes. By construction, $\overline\Phi$ factors into a product of $T_i^j$ for every $j$ for every $i$ and the boundary twists of $\pi_i$-holes for all $i$; we define $\Phi_i$ to be the product of all $T_i$-twists and the boundary twists of $\pi_i$-holes. We show that in the turned-over perspective the only hole which appears in more than one factor $\Phi_i$, is $\delta$. Hence, every other factorization splits into three factors, and these are completely determined by multiplicities.

Set $\pi^\text{in}$ as the outer boundary and consider the capped-off page $\mathbf D_\pi$ in the turned-over perspective. Here, no $\pi_i$-hole is encircled together with any $\pi_j$-hole for $i\neq j$, in symbols $m'_{\pi_i\pi_j}=0$, and $\delta$ is in at most $m'_{\delta}=k_1+k_2+2$ twists. Indeed, in the initial perspective (with $\nu^\text{out}$ as the outer boundary) the above multiplicities correspond to the fact that each twist which contains $\pi^\text{in}$ skips $\pi_i$-holes for one $i$ only, and that among the twists containing $\pi^\text{in}$ all $T_1$- and $T_2$-twists avoid $\delta$, which both holds by construction. 

On the other hand, the pairwise multiplicity of $\delta$ with $\pi^0_i$ in the turned-over perspective is exactly $m'_{\delta\pi^0_i}=k_i+1$. So, there are exactly $k_i+1$ twists encircling $\delta$ together with only $\pi_i$-holes. Again, this follows by construction as all the twists (in the initial perspective) which avoid any of the $\pi_i$-holes, avoid also $\pi^0_i$.

Therefore, since in the turned-over perspective there are no twists containing $\pi_i$ and $\pi_j$ together, we can consider the whole (abelianized) monodromy $\overline\phi$ as a product of three monodromies $\overline\phi_i$, uniquely determined by multiplicities:  the single and pairwise multiplicities of $\pi_i$-holes are the same as in $\overline\Phi$, and the twists around $\delta$ are distributed so that pairwise multiplicity of $\delta$ with the $\pi_i^0$-holes is $k_1+1,k_2+1$, and $0$, respectively. Thus, any positive factorization splits as $\overline{\Psi}_1\overline{\Psi}_2\overline{\Psi}_3$ with $\overline\Psi_i$ describing $\overline\phi_i$.
\endproof

For $\overline{\phi}_i\in\AbMap\mathbf D_{\pi}$ when $i=1$ or $2$, let us write out the given factorization $\overline{\Phi}_i$ as $F^0_i\cdots F^{l_i}_i$ and the arbitrary positive factorization $\overline{\Psi}_i$ as $P^0_i\cdots P^{l'_i}_i$. We order the Dehn twist factors so that viewed in $\mathbf D_{\pi}$ the ones containing $\pi^\text{in}$ come first and in the decreasing order of the number of holes they include. For $\overline{\Phi}_i$ the first $k_i+1$ twists $F^0_i,\dots,F^{k_i}_i$ then exactly equal $T^0_i,\dots,T^{k_i}_i$, and they are followed by the boundary twists of $\pi_i$-holes. For any positive factorization $\Psi_i$ we see that its length $l'_i+1$ is at least the length of the corresponding leg, $k_i+1$, and that $\pi^\text{in}$ is contained exactly in its first $k_i+1$ twists $P^0_i,\cdots, P^{k_i}_i$.

\begin{lem}\label{l1.2}
When the factorization $\overline{\Psi}_i$ does not agree with $\overline{\Phi}_i$, denote the first index on which they differ by $x_i:=\min\{j;F^j_i\neq P^j_i\}$. Then $x_i\leq k_i$, the holes in $\mathbf D_\pi$ encircled  by $P^{x_i}_i$ constitute a strict subset of the holes encircled by $F^{x_i}_i$, and no $\pi_i$-hole outside of $F^{x_i}_i$ is encircled by any non-boundary twist $P^j_i$ for $j\geq x_i$.
\end{lem}

\proof
Look at the capped-off page $\mathbf D_\pi$ in the turned-over perspective, with $\pi^\text{in}$ as the outer boundary. Containing $\pi^\text{in}$ in $\mathbf D_\pi$ exactly corresponds to containing $\delta$ in the turned-over perspective. Furthermore, for the twists containing $\pi^\text{in}$ in $\mathbf D_\pi$ the subset relation for encircled holes turns when viewed in the turned-over perspective. 

We first notice that $x_i$ always occurs among the twists containing $\delta$, thus $x_i\leq k_i$, as otherwise all pairwise multiplicities are reached within $\{P^j_i; 0\leq j\leq k_i\}$, and the factorization agrees with $\overline\Phi_i$. Now, if $P^{x_i}_i$ did not contain some hole $\chi$ contained in $F^{x_i}_i$, the pairwise multiplicity of $\chi$ with $\delta$ in $\overline\Psi_i$ would be strictly smaller than in $\overline\Phi_i$, in symbols $m'_{\chi\delta}(\overline\Psi_i)<m'_{\chi\delta}(\overline\Phi_i)$. Indeed, the number of twists containing $\delta$ is fixed and equal to $k_i+1$, and $F^j_i$ for $j\geq x_i$ all contain $\chi$, while $F^j_i$ for $j<x_i$ contains $\chi$ if and only if $P^j_i$ does. 

Finally, as in $\mathbf D$ the pairwise $M$-multiplicities of holes out of $F^{x_i}_i$ with any other $\pi$-hole are exactly as many as there are twists from $\{P^j_i=F^j_i;j<x_i\}$ around them, neither can be encircled together with any other $\pi$-hole additionally.
\endproof

\vspace{.2cm}\noindent\emph{Properties of the twists encircling $\nu$-holes.}\vspace{.1cm}

Having understood positive factorizations among the $\pi$-holes, the problem of finding a positive factorization reduces to whether any factorization $\overline{\Psi}$ (maybe $\overline\Phi$) of $\overline\phi\in\AbMap \mathbf D_\pi$ can be lifted to a positive factorization of $\phi\in\AbMap \mathbf D$.

In the following, we investigate possible lifts of the twists in $\overline\Psi$. In particular, we notice that the multiplicities put some restraints on the lifted twists, culminating in a list of properties satisfied by any positive factorization.

\begin{notation}
Let $n_i$ be the index of the first unknot on the leg $L_i$ whose stabilizations are not all of the same sign as the stabilizations of its leading unknot; when the leading unknot admits positive and negative stabilizations, we choose $n_i=0$, and when all stabilizations on $L_i$ are of the same sign, we set $n_i=k_i+1$.
\end{notation}

\begin{lem}\label{l2}
In any positive factorization of $\phi$ lifting $\overline\Psi$, there are at least $k_i-n_i+1$ twists among $\{P^{0}_i,\dots,P^{k_i}_i\}$ for $i=1$ and $2$, which lift to the twists which additionally contain only $\nu_i$-holes. They all contain $\pi^\text{in}$ and avoid $\pi^k_i$ for $k\leq n_i$; in particular, for $\overline\Phi_i$ these are exactly $\{F_i^{n_i},\dots,F_i^{k_i}\}$.
\end{lem}

\proof
Recall that on the disk with the initial outer boundary all the multiplicities $m_{\nu_i\nu_j}=0$ for $i\neq j$. In $\mathbf D$, this means that whenever some $\nu_i$ is encircled together with any of $\nu_j$, the twist needs to contain also the initial outer boundary, the hole $\nu^\text{out}$. On the other hand, for any $\nu_1$ or $\nu_2$ the pairwise multiplicity with $\pi^\text{in}$ in $\mathbf D$, $M_{\pi^\text{in}\nu_1}$ or $M_{\pi^\text{in}\nu_2}$, is greater than $M_{\pi^\text{in}\nu^\text{out}}=1$; precisely, for a $\nu^j_i$-hole the multiplicity is exactly $M_{\pi^\text{in}\nu^j_i}=k_i-j+2$. Indeed, in $\Phi$ there is a single twist around both $\pi^\text{in}$ and $\nu^\text{out}$ (which is the boundary twist of the outer boundary $\delta$), while $\nu^j_i$ are encircled together with $\pi^\text{in}$ also by the $T^k_i$-twists for $k\geq j$ (there are $k_i+1-j$ of them).
Thus, when there is some negative stabilization on $L_i$ (at least) $k_i-n_i+1$ Dehn twists which contain $\pi^\text{in}$ need to lift into twists which additionally include only the $\nu_i$-holes. 

Moreover, as on the initial disk the multiplicities $m_{\pi^k_i\nu_i}=0$ for $k\leq n_i$, whenever such $\pi^k_i$ is in $\mathbf D$ encircled together with $\nu_i$, the twist contains also $\nu^\text{out}$; hence, the $k_i-n_i+1$ twists mentioned above avoid all $\pi^k_i$ for $k\leq n_i$. In case of $\overline\Phi_i$ these are exactly $\{F_i^{n_i},\dots,F_i^{k_i}\}$.
\endproof

\begin{lem}\label{r2} 
In any positive factorization of $\phi$ lifting $\overline\Psi$, the lifts of twists from $\overline\Psi_3$ do not encircle any of $\nu^k_3$ for $k\leq n_3$. Furthermore, around $\nu^\text{out}$ and each $\nu^k_3$ for $k>n_3$, there are at least as many lifts of twists from $\overline\Psi_3$ as there are $T_3$-twists around it in the original factorization.
\end{lem}

\proof
Viewed in the turned-over perspective, the multiplicities $m'_{\nu^k_3\pi_3}=0$ for $k\leq n_3$. Hence, whenever such $\nu^k_3$ is in $\mathbf D$ encircled together with any $\pi_3$, the twist contains also $\pi^\text{in}$. But since no twist in $\overline\Psi_3$ contains $\pi^\text{in}$ in $\mathbf D$, their lifts necessarily avoid all $\nu^k_3$ for $k\leq n_3$.

Moreover, the pairwise multiplicity $M_{\nu^j_3\pi^k_3}=1+j-k$ for $j>k\geq n_3$, while $M_{\nu^j_3\pi^\text{in}}=1$. Hence, since all twists in $\overline\Psi_1$ and $\overline\Psi_2$ which contain any (hence all) $\pi_3$-hole contain also $\pi^\text{in}$, all but one twist which contain a $\nu_3$-hole together with some $\pi_3$, need to arise as lifts of twists in $\overline\Psi_3$; for $\nu^j_3$, there are $j-n_3$ of them, which is exactly the number of $T_3$-twists around $\nu^j_3$.
\endproof

\begin{prop}\label{pr2}
Properties of a positive factorization of $\phi$ concerning the holes $\nu^\text{out}\cup\nu_3$ in $\mathbf D$:
\begin{enumerate}[leftmargin=.6cm]
\item\label{p1} The pairwise multiplicity of each $\nu^\text{out}\cup\nu_3$ with any of $\pi^\text{in}\cup\pi_1\cup\pi_2$ is one, and with the $\pi_3$-holes the multiplicities equal $M_{\nu^j_3\pi^k_3}=\max\{1,1+j-k\}$.
\item\label{p2} Each of $\nu^j_3$-holes is encircled by at most $j+2$ non-boundary Dehn twists, $\nu^\text{out}$ by at most $k_3+2$.
\item\label{p3} The pairwise multiplicity of each $\nu^j_3$ with any $\nu^{\geq j}_3$ is exactly $j+1$.
\item\label{p4} Lemma \ref{l2}. Lifts of $k_i-n_i+1$ twists among $\{P^{0}_i,\dots,P^{k_i}_i\}$ for $i=1,2,$ never encircle any of $\nu^\text{out}\cup\nu_3$.
\item\label{p5} Lemma \ref{r2}. For every $\nu^\text{out}\cup\nu_3$ there is exactly one twist which encircles it together with any (and hence all) of $\pi_3$, and is not a lift of a twist from $\overline\Psi_3$. There is no lift of twists from $\overline\Psi_3$ around any $\nu^j_3$ with $j\leq n_3$, there are at least $j-n_3$ of them around $\nu^j_3$ with $j>n_3$, and at least $k_3-n_3+1$ around $\nu^\text{out}$.
\end{enumerate}
\end{prop}

\proof
Properties (\ref{p1})--(\ref{p3}) are obtained by counting twists in the original factorization $\Phi$. Property (\ref{p4}) restates Lemma \ref{l2}, and property (\ref{p5}) restates Lemma \ref{r2}.
\endproof

Taking the properties (\ref{p1}) and (\ref{p5}), we see that for each $\nu^\text{out}\cup\nu_3$ the twists in $\overline\Psi_1\cup\overline\Psi_2$ whose lifts encircle that hole, regarded as subsets of holes they encircle, form a partition of $\pi$-holes. Furthermore, according to the property (\ref{p2}) there is a bound on the number of parts -- that is, twists -- such a partition can consist of. Finally, the property (\ref{p3}) specifies how the partitions associated to different $\nu_3$-holes interact.

\begin{definition}
Let a {\em $\overline\Psi$-partition} be a set of twists from $\overline\Psi_1\cup\overline\Psi_2$ which -- as a subsets of holes they encircle -- partition $\pi$-holes. If two sets of twists from  $\overline\Psi_1\cup\overline\Psi_2$ define set-wise the same partition, then, since the twists of the two sets need to be parallel or equal, we say that the two $\overline\Psi$-partitions are {\em parallel}. The equal twists of different partitions are refered to as {\em shared}.

As noticed above, in a positive factorization every $\nu_3$-hole $\chi$ defines a $\overline\Psi$-partition; we will call it a {\em $\overline\Psi$-partition associated to $\chi$}.
\end{definition}

\vspace{.2cm}\noindent\emph{Lifting positive factorizations on $\pi$-holes over $\nu$-holes.}\vspace{.1cm}

We proceed successively, focusing on $\nu^j_3$ for $j$ in $0,1,\dots,k_3+1$; here, we denote $\nu^{k_3+1}_3:=\nu^\text{out}$. Let $\{j_l\}$ be a subsequence such that $|\nu^{j_l}_3|\neq 0$ (note that we have renamed the outer boundary of $\mathbf D$ to $\delta$, so it is not counted in $|\nu^0_3|$). Formally, we set 
\[\begin{array}{l}j_0=-1 \\
j_l:=\min\{j;j>j_{l-1}\text{ and }|\nu^j_3|\geq1\}. \end{array} \]

\begin{definition}
We will refer to the subsequence counter as the {\em level}. 
We will say that $\overline\Psi$ {\em lifts over $\nu_3^{\leq j_\ell}$} if the twists in $\overline\Psi$ can be lifted to  $\AbMap \mathbf D_{\pi\cup\nu^{\leq {j_\ell}}_3}$ so that the factorization satisfies the properties of Proposition \ref{pr2}.
\end{definition}

Note that to have a positive factorization, it is necessary that some $\overline\Psi$ lifts over all $\nu_3\cup\nu^\text{out}$.

\begin{lem}\label{l3.1}
If $\overline\Psi$ lifts over $\nu_3\cup\nu^\text{out}$, the $\overline\Psi$-partitions associated to every $\nu_3$-hole and to $\nu^\text{out}$ are all built of twists from the same  $\overline\Psi_i$ for either $i=1$ or $2$.
\end{lem}

\proof
As observed under Proposition \ref{pr2}, the set of all twists from $\overline\Psi_1\cup\overline\Psi_2$ whose lifts encircle a $\nu^\text{out}\cup\nu_3$-hole forms a partition of  the $\pi$-holes. 
Since every partition needs a twist which contains $\pi^\text{in}$, every $\overline\Psi$-partition consists of a twist $P^k_i$ with $k\leq k_i$ for either $i=1$ or $2$, and some twists covering all $\pi_i$-holes which are not encircled by $P^k_i$. Now, we separate three cases:
\begin{enumerate}[leftmargin=.7cm]
 \item[(i)] If all partitions have more than $j_1+2$ parts, the second (\ref{p2}) property can never be satisfied and there is no positive factorization. 
 
  \item[(ii)] If there is a partition of exactly $j_1+2$ parts, $j_1+1$ of them are necessarily shared by partitions associated to all $\nu_3\cup\nu^\text{out}$, to fulfill the third (\ref{p3}) property. Since around each $\pi$-hole there can be only one twist which does not contain $\pi^\text{in}$, the twists other than $P^k_i$ with $k\leq k_i$ are always shared by all $\nu_3\cup\nu^\text{out}$. Hence, the $\overline\Psi$-partitions associated to different $\nu^{j}_3$-holes come from different partitioning of holes contained in $P^k_i$, which is possible only by $\overline\Psi_i$-twists.
  
   \item[(iii)] If there is a partition of less than $j_1+2$ parts, we can extend all its defining twists over all $\nu^\text{out}\cup\nu_3$. Indeed, this choice satisfies the second (\ref{p2}) and the third (\ref{p3}) property (after we complete the factorization by lifts of the twists from $\overline\Psi_3$ and some twists which do not contain any $\pi$-holes), and the lifted twists obviously come from a single $\overline\Psi_i$. 
\endproof
\end{enumerate}

From the above proof we see that, in case (i) the factorization $\overline\Psi$ does not lift over $\nu_3^{j_1}$ and hence there is no positive factorization lifting $\overline\Psi$, while in case (iii) the factorization $\overline\Psi$ lifts over $\nu_3\cup\nu^\text{out}$ and there is no obstruction for positive factorization in terms of Proposition \ref{pr2}. In case (ii) the factorization $\overline\Psi$ lifts over $\nu_3^{j_1}$ and in this case, we continue by repeating a similar analysis for the holes $\nu_3^{j_l}$ with $l>1$. 

We first notice that, when looking for obstructions of positive factorization, it suffices to examine one particular factorization of $\overline\phi$, namely $\overline\Phi$.

\begin{notation}
We name the truncated continued fraction which corresponds to maximal fully positive (for $i=1,2$) or maximal fully negative (for $i=3$) truncation of the leg $L_i$, by $-\frac{1}{q_i}:=[a^0_i,\dots,a^{n_i-1}_i]$, or $-\frac{1}{q_i}:=-\infty$ when $n_i=0$.
\end{notation}

\begin{lem}\label{l3.2}
Assume that the legs are ordered so that $-\frac{1}{q_1}\geq -\frac{1}{q_2}$. If any positive factorization $\overline\Psi_i$ of either $\overline\phi_1$ or $\overline\phi_2$ lifts over $\nu^{\leq j}_3$, so does $\overline\Phi_1$.
\end{lem}

\proof
Suppose we are lifting $\overline\Psi_i$ which differs from $\overline\Phi_i$. At each level $l$ we are looking for partitions of the least possible parts among the partitions not associated to any $\nu_3^{<j_l}$. As observed in the proof of Lemma \ref{l3.1}, the partitions associated to $\nu_3^{j_l}$-holes are built from the partitions associated to $\nu_3^{j_{l-1}}$-holes by splitting the twist which contains $\pi^\text{in}$. This means that the twists $P^k_i$ for $k\leq k_i$ are being chosen for successive partitions in order of increasing $k$ (not necessary all of them). Since -- according to Lemma \ref{l1.2} -- the holes out of $F^{x_i}_i$ for $x_i=\min\{k;P^k_i\neq F^k_i\}$ are in any positive factorization $\overline\Psi_i$ encircled only by the twists which agree with some twists in $\overline\Phi_i$, the first $x_i$ of the $\overline\Psi_i$-partitions agree with the $\overline\Phi_i$-partitions; they are composed of a twist $P^k_i=F^k_i$ for $k< x_i$ and the boundary twists of the holes out of $P^k_i$, and the two factorizations lift simultaneously. But, once, at the level $\ell$, we use a $\overline\Psi_i$-partition which involves a twist $P^k_i\neq F^k_i$ with $k\leq k_i$, Lemma \ref{l1.2} tells that $\overline\Phi_i$ admits at least one more partition of at least one less part. Since by assumption $\overline\Psi_i$ lifts over $\nu^{\leq j_\ell}_3$, this $\overline\Phi_i$-partition has less than $j_\ell+2$ parts, and can be used for all $\nu^{\geq j_\ell}_3$, fulfilling the properties of Proposition \ref{pr2}.

According to Lemma \ref{l2}, every $\overline\Phi_i$-partition consists of a $T^k_i$-twist for $k<n_i$ and the boundary twists of the holes out of $T^k_i$. 
Comparing $\overline\Phi_1$ to $\overline\Phi_2$, the inequality $-\frac{1}{q_1}\geq -\frac{1}{q_2}$ means that at the first index in which the two continued fractions disagree, the coefficient $a^k_1$ is smaller than $a^k_2$. Hence, the corresponding $T^k_1$-twist avoid less $\pi$-holes than $T^k_2$ does, and the $(k+1)^\text{th}$ $\overline\Phi_1$-partition has less parts than the $(k+1)^\text{th}$ $\overline\Phi_2$-partition, while the first $k$ partitions have the same number of parts. 
\endproof

Lemma \ref{l3.2} essentially means that, when looking for obstructions of positive factorization, we can focus only on $\overline\Phi$ among $\overline\phi$-factorizations. Moreover, once we number the legs so that $-\frac{1}{q_1}\geq -\frac{1}{q_2}$, it suffices to check whether $\overline\Phi_1$ lifts over $\nu^{\leq j_l}_3$ for all levels $l$. If it does not, then no positive factorization $\overline\Psi_i$ of either $\overline\phi_1$ or $\overline\phi_2$ does (Lemma \ref{l3.2}). Hence, also $\overline\Psi$ cannot lift to a positive factorization of $\phi$ (Lemma \ref{l3.1}).

\begin{lem}\label{l4} 
Assume that $\overline\Phi_1$ lifts over $\nu_3^{<j_\ell}$ and that no $\overline\Phi_1$-partition associated to $\nu_3^{<j_\ell}$ could be extended over $\nu_3^{\geq j_\ell}$. 
At the $\ell^\text{th}$ level when $j_\ell< n_3$ (where we consider $n_3=k_3+1$ if there is no positive stabilization on $L_3$) one of the following happens:
\begin{enumerate}[leftmargin=.8cm]
  \item[(i)] If there is no $\overline\Phi_1$-partition into less than $j_\ell+2$ parts which has not been associated to some $\nu_3^{<j_\ell}$, and there are less than $|\nu^{j_\ell}_3|$ of parallel $\overline\Phi_1$-partitions into $j_\ell+2$ parts, there is no positive factorization of $\phi$.
  
  \item[(ii)] If there is no $\overline\Phi_1$-partition into less than $j_\ell+2$ parts which has not been associated to some $\nu_3^{<j_\ell}$, and there are at least $|\nu^{j_\ell}_3|$ of parallel $\overline\Phi_1$-partitions into $j_\ell+2$ parts, the factorization $\overline\Phi$ lifts over $\nu_3^{\leq j_\ell}$ and there are truncations of the legs $L_1$ and $L_3$ which are related as either: 
 \[\begin{array}{lllcccc}L_3^{(\ell)}&:&(-a_3^{j_1},&-2^{\times (j_2-j_1-1)},& -a_3^{j_2},&\dots,&-a_3^{j_\ell})\\ 
 L_1^{(\ell)}&:&(-3,-2^{\times (a_3^{j_1}-4)},&-j_2+j_1-2,&-2^{\times (a_3^{j_2}-3)},&\dots,& -2^{\times (a_3^{j_\ell}-3)})\end{array}\]
 or:
  \[\begin{array}{lllcccc}L_3^{(\ell)}&:&(-3,-2^{\times (j_1-1)},&-a_3^{j_1},&-2^{\times (j_2-j_1-1)},&\dots,&-a_3^{j_\ell})\\ 
 L_1^{(\ell)}&:&(-j_1-3,&-2^{\times (a_3^{j_1}-3)}&-j_2+j_1-2,&\dots,& -2^{\times (a_3^{j_\ell}-3)})\end{array}\]
  
  \item[(iii)] If there is a $\overline\Phi_1$-partition into less than $j_\ell+2$ parts, the factorization $\overline\Phi$ lifts over all $\nu_3\cup\nu^\text{out}$, but the surgery presentation contains a balanced sublink. 
\end{enumerate}
\end{lem}

\proof
The assumption means that $\overline\Phi_1$ falls under (ii) for all levels up to the $\ell^\text{th}$.  

At the $\ell^\text{th}$ level, if (i) there are only partitions of more than $j_\ell+2$ twists or there are less than $|\nu^{j_\ell}_3|$ of $j_\ell+2$-part partitions, there is no positive factorization; because we cannot satisfy the first (\ref{p1}) and the second (\ref{p2}) property of Proposition \ref{pr2} simultaneously. 

On the other hand, the conditions of (ii) allow us to obtain a positive factorization in $\mathbf D_{\pi\cup\nu_3^{\leq j_\ell}}$, but these conditions also prescribe how the truncations of legs are related. Indeed, $j_l$ always gives the index of an unknot on $L_3$ with surgery coefficient less than $-2$, $j_l-j_{l-1}$ counts the number of parallel twists, which is one more than the number of unknots with coefficient $-2$ preceding the unknot $v^{j_l}_3$; so, the corresponding part of $L_3$ looks like 
\[ \dots, -2^{\times(j_l-j_{l-1}-1)}, -a_3^{j_l},\dots \]
The fact that the conditions of (ii) are satisfied for the levels up to $\ell^\text{th}$ means that the number of separated holes in the partitions associated to $\nu_3^{j_l}$ compared to the partitions associated to $\nu_3^{j_{l-1}}$ is exactly $j_l-j_{l-1}$ (because there are no partitions of less than $j_l+2$ parts, there exists $j_l+2$-part partition, and previous partitions have $j_{l-1}+2$ parts), which on $L_1$ corresponds to an unknot of coefficient $-j_l+j_{l-1}-2$, which is followed by exactly $|\nu^{j_l}_3|-1$ of unknots with coefficient $-2$ (because there is at least $|\nu^{j_l}_3|$ partitions into $j_l+2$ parts at the $l^\text{th}$ level, and there is no $j_l+2$-part partition at the $(l+1)^\text{th}$ level). The corresponding part of $L_1$ looks like 
\[ \dots, -j_l+j_{l-1}-2, -2^{\times(|\nu^{j_l}_3|-1)}, \dots \]

Finally, the condition (iii) at the $\ell^\text{th}$ level quit the above sequence of pairing between coefficients of $L_3$ and $L_1$. We have $j_\ell-j_{\ell-1}$ parallel twists (so, $j_\ell-j_{\ell-1}-1$ of unknots with coefficient $-2$ on $L_3$) but we leave out less than $j_\ell-j_{\ell-1}$ holes by the next $T_{1}$-curve (its coefficient being at least $-j_\ell+j_{\ell-1}-1$). Since $j_\ell< n_3$ (hence, the preceding unknots of coefficient $-2$ being at the indices lower than $n_3$) and since according to Lemma \ref{l2} only the twists $T_1^0,\dots, T^{n_1-1}_1$ form $\overline\Phi_1$-partitions associated to any $\nu_3$ (hence, the unknot of coefficient greater than $-j_\ell+j_{\ell-1}-2$ being at the index lower than $n_1$), the two truncated chains correspond to the rational numbers smaller than or equal to $-\frac{1}{q_3}$ and $-\frac{1}{q_1}$, respectively. Comparing to Remark \ref{rmk:dual}, we see that $q_3$ and $q_1$ add up to at least $1$.
\endproof

\proof[Proof of Proposition \ref{prop:no+factor}]
 The process of lifting $\overline\Phi$ over $\nu_3$-holes eventually stops as we run into an obstruction for positive factorization (i) or we leave the assumed conditions (iii), if not before when we hit the $n_3$-level (the $(k_3+1)$-level if there is no positive stabilization on $L_3$). In the latter case, when (ii) is fulfilled by all levels $j_l< n_3$, we look at the possibilities of encircling $\nu^\text{out}$. 
 In order for a positive factorization to exist, we would -- according to (\ref{p2}) and (\ref{p5}) of Proposition \ref{pr2} -- need after the last checked level $\ell=\max\{l;j_l<n_3\}$ another $\overline\Phi$-partition of at most $n_3+1$ parts which has not been associated to any $\nu_3^{\leq j_\ell}$. But, the existence of such a partition would, as in the last paragraph of the previous proof, imply that there is a balanced sublink in the surgery presentation. Indeed, there would be $n_3-j_\ell-1$ of unknots with coefficient $-2$ preceding $v^{n_3}_3$ on $L_3$, and the corresponding unknot on the truncated $L_1$ would have coefficient at least $-n_3+j_\ell-1$, meaning by Remark \ref{rmk:dual} that $q_3+q_1\geq 1$.
\endproof

\proof[Proof of Theorem \ref{thm} and Proposition \ref{prop:tightSxS}]
Joining Proposition \ref{prop:+factor} and Proposition  \ref{prop:no+factor} we obtain the theorem, and the proposition should be read as its special case. 

Indeed, Legendrian surgeries on the tight $S^1\times S^2$, given by a balanced link (as in Proposition \ref{prop:+factor}), give Stein fillable structures. On the other hand, in the absence of the balanced sublink, Proposition \ref{prop:no+factor} tells that the associated planar monodromy do not admit positive factorization and hence, because of the Wendl's theorem (Theorem \ref{thm:W}), the presented contact manifold do not admit any Stein filling. 
\endproof

%%%%%%%%%%%%%%%%%%%%%%%%%%%%%%%%%%%%%%%%%%%%%%%%%%%%%%%%%%%%%%

%%%%%%%%%%%%%%%%%%%%%%%%%%%%%%%%%%%%%%%%%%%%%%%%%%%%%%%%%%%%%%
%%%%%%%%%%%%%%%%%%%%%%%%%%%%%%%%%%%%%%%%%%%%%%%%%%%%%%%%%%%%%%

\begin{thebibliography}{99}

\bibitem{DG}
F. Ding, and H. Geiges, {\em Handle moves in contact surgery diagrams}, J. Topol. \textbf{2} (2009) 105--122.

\bibitem{EMVH-M}
H. Endo, T. Mark, and J. Van Horn-Morris, {\em Monodromy substitutions and rational blowdowns}, J. Topol. \textbf{4} (2011),  227--253.

\bibitem{G}
P. Ghiggini, {\em On tight contact structures with negative maximal twisting number on small Seifert manifolds}, Algebr. Geom. Topol. \textbf{8} (2008) 381--396.

\bibitem{GLS0}
P. Ghiggini, P. Lisca, and A. Stipsicz, {\em Classification of tight contact structures on small Seifert 3-manifolds with $e_0\geq 0$}, Proc. Amer. Math. Soc. \textbf{134} (2006) 909--916.

\bibitem{GLS}
P. Ghiggini, P. Lisca, and A. Stipsicz, {\em Tight contact structures on some small Seifert fibered 3-manifolds}, Amer. J. Math. \textbf{129}(\textbf{5}) (2007) 1403--1447.

\bibitem{Gompf}
R. Gompf, {\em Handlebody construction of Stein surfaces}, Ann. of Math. \textbf{148} (1998) 619--693.

\bibitem{LL}
A. G. Lecuona, and P. Lisca, {\em Stein fillable Seifert fibered 3–manifolds}, Algebr. Geom. Topol. \textbf{11} (2011) 625--642.

\bibitem{LS.III}
P. Lisca, and A. Stipsicz, {\em Ozsváth-Szabó invariants and tight contact 3-manifolds III}, J. Symplectic Geom. \textbf{5}(\textbf{4}) (2007) 357--384.

\bibitem{MMc}
D. Margalit, and J. McCammond, {\em Geometric presentations for the pure braid group}, J. Knot Theory Ramifications \textbf{18} (2009) 1--20.

\bibitem{M}
I. Matkovi\v{c}, {\em Classification of tight contact structures on small Seifert fibered $L$-spaces}, Algebr. Geom. Topol. \textbf{18} (2018) 111--152.

\bibitem{PVH-M}
O. Plamenevskaya, and J. Van Horn-Morris, {\em Planar open books, monodromy factorizations and symplectic fillings}, Geom. Topol. \textbf{14} (2010) 2077--2101.

\bibitem{Sch}
S. Schönenberger, {\em Determining symplectic fillings from planar open books}, J. Symplectic Geom., \textbf{5}(\textbf{1}) (2007) 19--41.

\bibitem{S}
A. Stipsicz, {\em Ozsváth-Szabó invariants and 3-dimensional contact topology}, Proceedings of the International Congress of Mathematicians, Hyderabad 2010, vol II, 1159--1178.

\bibitem{W}
C. Wendl, {\em Strongly fillable contact manifolds and $J$-holomorphic foliations}, Duke Math. J. \textbf{151} (2010) 337--384.

\bibitem{Wu}
H. Wu, {\em Legendrian vertical circles in small Seifert spaces}, Commun. Contemp. Math. \textbf{8} (2006) 219--246.

\end{thebibliography}
\end{document}